\documentclass[12pt,vatola]{article}

\usepackage{graphicx}

\def\C{\centerline}

\def\re#1{\par\hangindent\parindent\indent\llap{#1\enspace}\ignorespaces}

\def\no{\noindent}

\begin{document}

\C{\large\bf The Number of Complete Maps on Surfaces} \vskip 5mm

\C{\Large\bf }

\vskip 10mm \C{Linfan Mao} \vskip 3mm \C{\scriptsize (Chinese
Academy of Mathematics and System Science, Beijing 100080,
P.R.China)}

\vskip 5mm \C{Yanpei Liu} \C{\scriptsize (Institute of Applied
Mathematics, Northern Jiaotong University, Beijing,100044) }

\vskip 5mm \C{Feng Tian}\vskip 3mm \C{\scriptsize (Chinese Academy
of Mathematics and System Science, Beijing 100080, P.R.China)}

\vskip 8mm
\begin{minipage}{130mm}
\no{{\bf Abstract}: {\small A map is a connected topological graph
cellularly embedded in a surface and a complete map is a
cellularly embedded complete graph in a surface. In this paper,
all automorphisms of complete maps of order $n$ are determined by
permutations on its vertices. Applying a scheme for enumerating
maps on surfaces with a given underlying graph, the numbers of
unrooted complete maps on orientable or non-orientable surfaces
are obtained. } }

\no{{\bf Key words:} {\small  embedding, complete map,
isomorphism, automorphism group, Burnside Lemma.}}

\no{{\bf Classification:}  AMS(2000) 05C10,05C25, 05C30}
\end{minipage}

\vskip 10mm

\no{\bf 1. Introduction} \vskip 5mm

\no All surfaces considered in this paper are 2-dimensional
compact closed manifolds without boundary, graphs are connected
and simple graphs with the maximum valency $\geq 3$ and groups are
finite. For terminologies and notations not defined here can be
seen in $[21]$ for maps, $[20]$ for graphs and in $[2]$ for
permutation groups.

The enumeration of rooted maps on surfaces, especially, the
sphere, has been intensively investigated by many researchers
after the Tutte's pioneer work in 1962 (see [$21$]). Comparing
with rooted maps, observation for the enumeration of unrooted maps
on surface is not much. By applying the automorphisms of the
sphere, Liskovets gave an enumerative scheme for unrooted planar
maps(see [$12$]). Liskovets, Walsh and Liskovets got many
enumeration results for {\it general planar maps, regular planar
maps, Eulerian planar maps, self-dual planar maps} and {\it
2-connected planar maps}, etc (see $[12]-[14]$).

General results for the enumeration of unrooted maps on surface
other than sphere are very few.  Using the well known Burnside
Lemma in permutation group theory, Biggs and White presented a
formula for enumerating non-equivalent embeddings of a given graph
on orientable surfaces$^{[2]}$, which are the classification of
embeddings by orientation-preserving automorphisms of orientable
surfaces. Following their idea, the numbers of non-equivalent
embeddings of complete graphs,complete bipartite graphs, wheels
and graphs whose automorphism group action on its ordered pair of
adjacent vertices is semi-regular are gotten in references
$[15]-[16],[20]$ and $[11]$. Although this formula is not very
efficient and need more clarifying for the actual enumeration of
non-equivalent embeddings of a graph, the same idea is more
practical for enumerating rooted maps on orientable or
non-orientable surfaces with given underlying graphs(see
$[8]-[10]$).

For projective maps with a given 3-connected underlying graph,
Negami got an enumeration result for non-equivalent embeddings by
establishing the double planar covering of projective maps(see
[$18$]). In [$7$], Jin Ho Kwak and Jaeun Lee obtained the number
of non-congruent embeddings of a graph, which is also related to
the topic discussed in this paper.

Combining the idea of Biggs and White for non-equivalent
embeddings of a graph on orientable surfaces and the Tutte's
algebraic representation for maps on surface$^{[19],[21]}$, a
general scheme for enumerating unrooted maps on locally orientable
surfaces with a given underlying graph is obtained in this paper.
Whence, the enumeration of unrooted maps on surfaces can be
carried out by the following programming:

\vskip 3mm
 STEP 1. Determined all automorphisms of maps with
a given underlying graph;

\vskip 2mm
 STEP 2. Calculation the the fixing set
$Fix(\varsigma)$ for each automorphism $\varsigma$ of maps;

\vskip 2mm
 STEP 3. Enumerating the unrooted maps on surfaces
with a given underlying graph by this scheme.

\vskip 2mm

Notice that this programming can be used for orientable or
non-orientable surfaces, respectively and get the numbers of
orientable or non-orientable unrooted maps underlying a given
graph.

The main purpose of this paper is to  enumerate the orientable or
non-orientable complete maps. In $1971$, Biggs proved$^{[1]}$ that
the order of automorphism group of an orientable complete map of
order $n$ divides $n(n-1)$, and equal $n(n-1)$ only if the
automorphism group of the complete map is a {\it Frobenius group}.
In this paper, we get a representation by the permutation on its
vertices for the automorphisms of orientable or non-orientable
complete maps. Then as soon as we completely calculate the fixing
set $Fix(\varsigma)$ for each automorphism $\varsigma$ of complete
maps, the enumeration of unrooted orientable or non-orientable
complete maps can be well done by our programming.

The problem of determining which automorphism of a graph is an
automorphism of a map is also interesting for {\it Riemann}
surfaces or {\it Klein} surfaces - surfaces equipped with an
analytic or dianalytic structure, for example, automorphisms of
Riemann or Klein surfaces have be given more attention since
1960s, see for example,$[3]-[4],[6],[17],$ but it is difficult to
get a concrete representation for an automorphism of Riemann or
Klein surfaces. The approach used in this paper can be also used
for combinatorial discussion automorphisms of Riemann or Klein
surface.

Terminologies and notations used in this paper are standard. Some
of them are mentioned in the following.

For a given connected graph $\Gamma$, an {\it embedding} of
$\Gamma$ is a pair $({\cal J},\lambda)$, where ${\cal J}$  is a
rotation system of $\Gamma$, and $\lambda : E(\Gamma) \rightarrow
Z_2$. The edge with $\lambda (e)= 0$ or $\lambda (e)=1$ is called
the {\it type $0$} or {\it type $1$ edge}, respectively.

A {\it{map}} $M = ({\cal X} _{\alpha,\beta},\cal{P})$ is defined
to be a permutation $\cal{P}$ acting on ${\cal X} _{\alpha,\beta}$
of a disjoint union of quadricells $Kx$ of $x\in X$, where, $X$ is
a finite set and $K=\{1,\alpha,\beta,\alpha\beta \}$ is the {\it
Klein} group, satisfying the following conditions:

$i)$ $\forall{x}\in {{\cal X} _{\alpha,\beta}}$, there does not
exist an integer $k$ such that ${\cal{P}}^{k}x = \alpha x$;

$ii)$ $\alpha{\cal{P}}={\cal{P}}^{-1}\alpha$;

$iii)$ the group $\Psi_{J}=\left<\alpha,\beta,\cal{P}\right>$ is
transitive on ${\cal X}_{\alpha,\beta}$.

According to the condition $ii)$, the vertices of a map are
defined to be the pairs of conjugate of ${\mathcal P}$ action on
${\mathcal X}_{\alpha,\beta}$ and edges the orbits of $K$ on
${\mathcal X}_{\alpha,\beta}$, for example, for $\forall x\in
{\mathcal X}_{\alpha,\beta}$, $\{x,\alpha x,\beta x,\alpha\beta
x\}$ is an edge of the map $M$. Geometrically, any map $M$ is an
embedding of a graph $\Gamma$ on a surface, denoted by
$M=M(\Gamma)$ and $\Gamma=\Gamma(M)$ ( see also $[19]-[21]$ for
details). The graph $\Gamma$ is called the underlying graph of the
map $M$. If $r\in {\mathcal X}_{\alpha,\beta}$ is marked
beforehand, then $M$ is called a {\it{rooted map}}, denoted by
$M^{r}$. A map is said non-orientable or orientable if the group
$\Psi_I=\left<\alpha\beta, {\mathcal P}\right>$ is transitive on
${\mathcal X}_{\alpha,\beta}$ or not.

For example, the graph $K_4$ on the tours with one face length $4$
and another $8$ , shown in the following Fig.$1$,

\includegraphics[bb=-10 10 200 230]{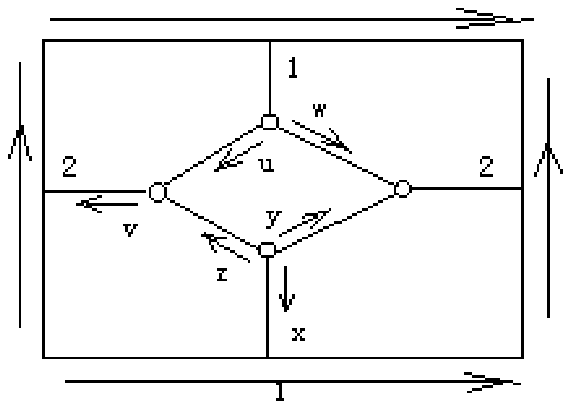}

\C{\bf Fig.$1$} \vskip 3mm \no{can be algebraically represented as
follows:}

{\it A map $({\mathcal X}_{\alpha,\beta},\mathcal{P})$ with
${\mathcal X}_{\alpha,\beta}= \{x,y,z,u,v,w,\alpha x,\alpha y,$ $
\alpha z,\alpha u,\alpha v,\alpha w, \beta x,\beta y,\beta z,\beta
u,$ $
 \beta v,\beta w,\alpha\beta x,
\alpha \beta y,\alpha \beta z,\alpha \beta u,\alpha \beta
v,\alpha\beta w \}$ and}

\begin{eqnarray*}
{\mathcal P} &=& (x,y,z)(\alpha \beta x,u,w)(\alpha \beta z,\alpha
\beta u,v)
(\alpha \beta y,\alpha \beta v,\alpha \beta w)\\
&\times& (\alpha x,\alpha z,\alpha y)(\beta x,\alpha w,\alpha u)(\beta z,\alpha v,\beta u)(\beta y,\beta w,\beta v)
\end{eqnarray*}

The four vertices of this map are $ \{(x,y,z), (\alpha x,\alpha
z,\alpha y)\}$, $\{(\alpha \beta x,u,w),(\beta x,\alpha w,\alpha
u)\}$, $\{(\alpha \beta z,\alpha \beta u,v),(\beta z,\alpha
v,\beta u)\}$ and $\{(\alpha \beta y,\alpha \beta v,\alpha \beta
w),(\beta y,\beta w,\beta v)\}$ and six edges are $\{e,\alpha
e,\beta e,\alpha\beta e\}$ for $\forall e\in \{x,y,z,u,v,w\}$.

Two maps $M_{1} = ({\mathcal
X}_{\alpha,\beta}^{1},{\mathcal{P}}_{1})$ and $M_{2} = ({\mathcal
X}_{\alpha,\beta}^{2},{\mathcal{P}}_{2})$ are said to be
{\it{isomorphic}} if there exists a bijection $\tau: {\mathcal
X}_{\alpha,\beta}^{1} \longrightarrow {\mathcal
X}_{\alpha,\beta}^{2}$ such that for $\forall{x}\in{{\mathcal
X}_{\alpha,\beta}^{1}}$,$\tau\alpha(x)=\alpha\tau(x)$,
$\tau\beta(x)=\beta\tau(x)$ and
$\tau{\mathcal{P}}_{1}(x)={\mathcal{P}}_{2}\tau(x)$. $\tau$ is
called an {\it{isomorphism}} between them. If $M_{1}=M_{2}=M$,
then an isomorphism between $M_{1}$ and $M_{2}$ is called an {\it
automorphism} of $M$. All automorphisms of a map $M$ form a group,
called the {\it{automorphism group}} of $M$ and denoted by ${\rm
AutM}$. Similarly, two rooted maps  $M^{r}_{1},$ $M^{r}_{2}$ are
said to be {\it{isomorphic}} if there is an isomorphism $\theta$
between them such that $\theta(r_{1})= r_{2}$, where $r_{1}$,
$r_{2}$ are the roots of $M_{1}^{r}$ , $M_{2}^{r}$, respectively
and denote the {\it{automorphism group}} of $M^{r}$ by ${\rm
AutM}^{r}$. It has been known that ${\rm AutM}^{r}$ is a trivial
group.

According to their action, isomorphisms between maps can divided
into two classes: {\it cyclic order-preserving isomorphism} and
{\it cyclic order-reversing isomorphism}, defined as follows,
which is useful for determining automorphisms of a map underlying
a graph.

For two maps $M_1$ and $M_2$, a bijection $\xi$ between $M_1$ and
$M_2$ is said to be {\it cyclic order-preserving} if for
$\forall{x}\in{{\mathcal
X}_{\alpha,\beta}^{1}}$,$\tau\alpha(x)=\alpha\tau(x)$,
$\tau\beta(x)=\beta\tau(x)$, $\tau{\mathcal{P}}_{1}(x)
={\mathcal{P}}_{2}\tau(x)$ and {\it cyclic order-reversing} if
$\tau\alpha(x)=\alpha\tau(x)$, $\tau\beta(x)=\beta\tau(x)$
$\tau{\mathcal{P}}_{1}(x)={\mathcal{P}}_{2}^{-1}\tau(x)$.

Now let $\Gamma$ be a connected graph. The notations ${\mathcal E
}^{O}(\Gamma ), {\mathcal E}^{N}(\Gamma )$ and ${\mathcal
E}^{L}(\Gamma )$ denote the embeddings of $\Gamma$ on the
orientable surfaces, non-orientable surfaces and locally surfaces,
${\mathcal M}(\Gamma )$ and ${\rm Aut}\Gamma$ denote the set of
non-isomorphic maps underlying $\Gamma$ and its automorphism
group, respectively.

\vskip 10mm

\no{\bf 2. The enumerative scheme for maps underlying a graph}
\vskip 6mm

\no A permutation $p$ on set $\Omega$ is called semi-regular if
all of its orbits have the same length. For a given connected
graph $\Gamma$, $\forall g\in{Aut \Gamma}$, $M=({\mathcal
X}_{\alpha,\beta},{\mathcal P})\in {\mathcal M}(\Gamma)$, define
an extended action of $g$ on $M$ to be

$$
g^*:{\mathcal X}_{\alpha,\beta} \longrightarrow {\mathcal
X}_{\alpha,\beta},
$$

\no{such that $M^{g^*}=gMg^{-1}$ with $g\alpha =\alpha g$ and
$g\beta = \beta g$. }

We have already known the following two results.

\vskip 4mm \no{\bf Lemma $2.1^{[21]}$}  {\it  For any rooted map
$M^r$, $AutM^r$ is trivial.}

\vskip 3mm \no{\bf Lemma $2.2^{[2],[21]}$} {\it For a given map
$M$, $\forall\xi\in AutM$, $\xi$ transforms vertices to vertices,
edges to edges and faces to faces on a map $M$, i.e, $\xi$ can be
naturally extended to an automorphism of surfaces. } \vskip 3mm

\vskip 4mm \no{\bf Lemma $2.3$} {\it If there is an isomorphism
$\xi$ between maps $M_1$ and $M_2$, then $\Gamma (M_1)=\Gamma
(M_2) =\Gamma$ and $\xi\in {\rm Aut\Gamma }$ if $\xi$ is cyclic
order-preserving  or $\xi\alpha\in {\rm Aut \Gamma }$ if $\xi$ is
cyclic order-reversing.} \vskip 3mm

{\it Proof} \ By the definition of an isomorphism between maps,
if $M_1=({\mathcal X}_{\alpha,\beta}^1,{\mathcal P}_1)$ is
isomorphic with $M_2=({\mathcal X}_{\alpha,\beta}^2,{\mathcal
P}_2)$, then there is an $1-1$ mapping $\xi$ between ${\mathcal
X}_{\alpha,\beta}^1$ and ${\mathcal X}_{\alpha,\beta}^2$ such that
$({\mathcal P}_1)^{\xi} = {\mathcal P}_2$ . Since isomorphic
graphs are considered to be equal, we get that $\Gamma
(M_1)=\Gamma (M_2)= \Gamma$. Now since

$$({\mathcal P}_2)^{-1} = ({\mathcal P}_2)^{\alpha}.$$

\no{We get that $\Gamma^{\xi}=\Gamma$ or $\Gamma^{\xi\alpha} =
\Gamma$, whence, $\xi\in {\rm Aut\Gamma}$  or $\xi\alpha\in {\rm
Aut\Gamma}$.\quad\quad $\natural$

According to Lemma 2.3, For $\forall g \in {\rm Aut}\Gamma,
\forall M \in {{\mathcal E}^{L}(\Gamma )}$, the induced action
$g^*$ of $g$ on $M$ is defined by $M^{g^*} =gMg^{-1} = ({\mathcal
X} _{\alpha,\beta}, g{\mathcal P}g^{-1})$.

Since ${\mathcal P}$ is a permutation on the set ${\mathcal X}
_{\alpha,\beta}$, by a simple result in permutation group theory,
${\mathcal P}^g$ is just the permutation replaced each element $x$
in $\mathcal P$ by $g(x)$. Whence $M$ and $M^{g^*}$ are
isomorphic. Therefore, we get the following enumerative theorem
for unrooted maps underlying a graph.

\vskip 4mm \no{{\bf Theorem 2.1} {\it For a connected graph
$\Gamma$, let ${\mathcal E} \subset {\mathcal E}^{L}(\Gamma )$.
Then the number $n({\mathcal E},\Gamma)$ of unrooted maps in
${\mathcal E}$ is}}

$$
n({\mathcal E},\Gamma )=\frac{1}{|{\rm Aut}\Gamma\times
\left<\alpha\right>|}\sum\limits_{g \in {\rm Aut}\Gamma\times
\left<\alpha\right>} |\Phi (g)|,
$$
\no{\it where, $\Phi(g)=\{ {\mathcal{P}}| {\mathcal{P}}\in
{\mathcal E}$ and ${\mathcal{P}}^g= {\mathcal{P}} \}$.} \vskip 3mm

 {\it Proof} \ According to Lemma 2.1, two maps $M_1,M_2\in {\mathcal E}$ are
 isomorphic if and only
if there exists an isomorphism $\theta\in{\rm Aut}\Gamma\times
<\alpha>$ such that $M_1^{\theta^*}=M_2$. Whence, we get that all
the unrooted maps in ${\mathcal E}$ are just the representations
of orbits in ${\mathcal E}$ under the action of ${\rm
Aut}\Gamma\times \left<\alpha\right>$. By the Burnside Lemma, we
get the following result for the number of unrooted maps in
${\mathcal E}$

$$
n({\mathcal E},\Gamma )=\frac{1}{|{\rm Aut}\Gamma\times
\left<\alpha\right>|}\sum\limits_{g \in {\rm Aut}\Gamma\times
\left<\alpha\right>} |\Phi (g)|.\quad\quad\natural
$$

\vskip 4mm \no{{\bf Corollary 2.1} {\it For a given graph
$\Gamma$, the numbers of unrooted maps in ${\mathcal E
}^{O}(\Gamma ),{\mathcal E}^{N}(\Gamma )$ and ${\mathcal
E}^{L}(\Gamma )$ are}}

$$
n^{O}(\Gamma )=\frac{1}{|{\rm Aut}\Gamma\times
\left<\alpha\right>|}\sum\limits_{g \in {\rm Aut}\Gamma\times
\left<\alpha\right>} |\Phi^{O} (g)|; \quad\quad (2.1)
$$

$$
n^{N}(\Gamma )=\frac{1}{|{\rm Aut}\Gamma\times
\left<\alpha\right>|}\sum\limits_{g \in {\rm Aut}\Gamma\times
\left<\alpha\right>} |\Phi^{N} (g)|; \quad\quad (2.2)
$$

$$
n^{L}(\Gamma )=\frac{1}{|{\rm Aut}\Gamma\times
\left<\alpha\right>|}\sum\limits_{g \in {\rm Aut}\Gamma\times
\left<\alpha\right>} |\Phi^{L} (g)|, \quad\quad (2.3)
$$

\noindent{\it where, $\Phi^{O}(g)=\{ {\mathcal{P}}|
{\mathcal{P}}\in {\mathcal E}^{O}(\Gamma )$ and ${\mathcal{P}}^g=
{\mathcal{P}} \}$, $\Phi^{N}(g)=\{ {\mathcal{P}}| {\mathcal{P}}\in
{\mathcal E}^{N}(\Gamma )$ and ${\mathcal{P}}^g= {\mathcal{P}}
\}$, $\Phi^{L}(g)=\{ {\mathcal{P}}| {\mathcal{P}}\in {\mathcal
E}^{L}(\Gamma )$ and ${\mathcal{P}}^g= {\mathcal{P}} \}$.} \vskip
4mm

\no{{\bf Corollary 2.2} {\it Let ${\mathcal E}(S,\Gamma )$ be the
embeddings of $\Gamma$ in the surface $S$,
 then the number $n(\Gamma ,S)$ of unrooted maps on $S$ with underlying g$\Gamma$ is }}
$$
n(\Gamma ,S)=\frac{1}{|{\rm Aut}\Gamma\times
\left<\alpha\right>|}\sum\limits_{g \in {\rm Aut}\Gamma\times
\left<\alpha\right>} |\Phi (g)|,
$$
\noindent{\it where, $\Phi(g)=\{ {\mathcal{P}}| {\mathcal{P}}\in
{\mathcal E}(S,\Gamma )$ and ${\mathcal{P}}^g= {\mathcal{P}} \}$.
}

 \vskip 3mm \no{{\bf
Corollary 2.3} {\it In formulae (2.1)-(2.3), $|\Phi (g)|\not=0$ i
and only if $g$ is an automorphism of an orientable or
non-orientable map underlying $\Gamma$.}}

\vskip 2mm

Directly using these formulae (2.1)-(2.3) to count unrooted maps
with a given underlying graph is not straightforward. More
observation should be considered. The following two lemmas give
necessary conditions for an induced automorphism of a graph
$\Gamma$ to be an cyclic order-preserving automorphism of a
surface.

\vskip 4mm \no{{\bf Lemma $2.4$} {\it For a map $M$ underlying a
graph $\Gamma$, $\forall g\in{\rm AutM}, \forall x\in{{\mathcal
X}_{\alpha,\beta}}$ with $X=E(\Gamma)$, }}

($i$)\quad $|x^{\rm AutM}|=|\rm AutM|$ ;

($ii$)\quad $|x^{< g >}|= o(g)$,

\no{where, $o(g)$ denotes the order of $g$.} \vskip 3mm

{\it Proof} \ For a subgroup $H < {\rm AutM}$, we know that
$|H|=|x^H||H_x|$. Since $H_x < {\rm AutM}^x$, where $M^x$ is a
rooted map with root $x$, we know that $|H_x|=1$ by Lemma $2.1$.
Whence, $|x^H|=|H|$. Now take $H={\rm AutM}$ or $< g
>$, we get the assertions ($i$) and ($ii$).\quad\quad $\natural$

\vskip 4mm \no{{\bf Lemma $2.5$} {\it Let $\Gamma$ be a connected
graph and $g\in{\rm Aut \Gamma}$. If there is a map $M\in
{\mathcal E}^{L}(\Gamma)$ such that the induced action $g^*
\in{\rm AutM}$, then for $\forall (u,v),(x,y)\in{E(\Gamma)}$,}}

$$[l^g(u),l^g(v)]=[l^g(x),l^g(y)]= constant,$$

\no{\it where, $l^g(w)$ denotes the length of the cycle containing
the vertex $w$ in the cycle decomposition of $g$ and $[a,b]$ the
least common multiple of integers $a$ and $b$.} \vskip 3mm

 {\it Proof} \ According to Lemma $2.4$, we know that the length of any
 quadricell $u^{v+}$ or
$u^{v-}$ under the action of $g^*$ is $[l^g(u),l^g(v)]$. Since
$g^*$ is an automorphism of map, therefore, $g^*$ is semi-regular.
Whence, we get that
 $$ [l^g(u),l^g(v)]=[l^g(x),l^g(y)]= constant.\quad\quad\natural$$

Now we consider conditions for an induced automorphism of a map by
an automorphism of graph to be a cyclic order-reversing
automorphism of surfaces.

\vskip 4mm \no{\bf Lemma $2.6$} {\it If $\xi\alpha$ is an
automorphism of a map, then $\xi\alpha =\alpha\xi .$} \vskip 3mm
{\it Proof} \ Since $\xi\alpha$ is an automorphism of a map, we
know that

$$(\xi\alpha)\alpha = \alpha (\xi\alpha).$$

\no{That is, $\xi\alpha =\alpha\xi .\quad\quad\natural$}

\vskip 4mm \no{\bf Lemma $2.7$} {\it If $\xi$ is an automorphism
of $M=({\mathcal X}_{\alpha,\beta},{\mathcal P})$, then $\xi\alpha
$ is semi-regular on ${\mathcal X}_{\alpha,\beta}$ with order
$o(\xi)$ if $o(\xi)\equiv 0(mod 2)$ or $2o(\xi)$ if $o(\xi)\equiv
1(mod 2)$.}

\vskip 3mm

{\it Proof} \ Since $\xi$ is an automorphism of map by Lemma
$2.6$, we know that the cycle decomposition of $\xi$ can be
represented by

$$\xi=\prod_k (x_1,x_2,\cdots,x_k)(\alpha x_1,\alpha x_2,\cdots, \alpha x_k),$$

\no{where, $\prod_k$ denotes the product of disjoint cycles with
length $k=o(\xi)$.}

Therefore,  if $k\equiv 0 (mod 2)$, we get that

$$\xi\alpha = \prod_k (x_1,\alpha x_2,x_3,\cdots,\alpha x_k)$$

\no{and if $k\equiv 1 (mod 2)$, we get that}

$$\xi\alpha = \prod_{2k} (x_1,\alpha x_2,x_3,
\cdots, x_k,\alpha x_1,x_2,\alpha x_3,\cdots,\alpha x_k).$$

\no{Whence, $\xi$ is semi-regular acting on ${\mathcal
X}_{\alpha,\beta}$.\quad\quad $\natural$}

\vskip 2mm

Now we can prove the following result for cyclic order-reversing
automorphisms of maps.

\vskip 4mm \no{\bf Lemma $2.8$}  {\it For a connected graph
$\Gamma$, let ${\mathcal K}$  be all automorphisms in ${\rm
Aut}\Gamma$ whose extending action on ${\mathcal
X}_{\alpha,\beta}$, $ X=E(\Gamma)$, are automorphisms of maps
underlying the graph $\Gamma$. Then for $\forall\xi\in {\mathcal
K}$, $o(\xi^*)\geq 2$, $\xi^*\alpha \in {\mathcal K}$ if and only
if $o(\xi^*)\equiv 0(mod 2).$}

 \vskip 3mm
{\it Proof} \ Notice that by Lemma $2.7$, if $\xi^*$ is an
automorphism of a map underlying $\Gamma$, then $\xi^*\alpha$ is
semi-regular acting on ${\mathcal X}_{\alpha,\beta}$.

Assume $\xi^*$ is an automorphism of the map $M=({\mathcal
X}_{\alpha,\beta},{\mathcal P})$. Without loss of generality, we
assume that

$$
{\mathcal P}=C_1C_2\cdots C_k,
$$
\no{where,$C_i =(x_{i1},x_{i2},\cdots,x_{ij_i})$ is a cycle in the
decomposition of $\xi |_{V(\Gamma)}$ and
$x_{it}=\{(e^{i1},e^{i2},\cdots,e^{it_i})(\alpha e^{i1},\alpha
e^{it_i},\cdots,\alpha e^{i2})\}$,}

$$\xi |_{E(\Gamma)}=(e_{11},e_{12},\cdots,e_{s_1})(e_{21},e_{22},\cdots,e_{2s_2})\cdots
(e_{l1},e_{l2},\cdots,e_{ls_l}),$$

\no{and}
$$
\xi^*=C (\alpha C^{-1}\alpha ),
$$

\no{where,
$C=(e_{11},e_{12},\cdots,e_{s_1})(e_{21},e_{22},\cdots,e_{2s_2})\cdots
(e_{l1},e_{l2},\cdots,e_{ls_l})$. Now since $\xi^*$ is an
automorphism of a map, we get that $s_1=s_2=\cdots =s_l
=o(\xi^*)=s.$}

If $o(\xi^*)\equiv 0(mod 2)$, define a map $M^*=({\mathcal
X}_{\alpha,\beta},{\mathcal P}^*)$ with

$$
{\mathcal P}^*=C_1^*C_2^*\cdots C_k^*,
$$
\no{where, $C_i^* =(x_{i1}^*,x_{i2}^*,\cdots,x_{ij_i}^*)$,
$x_{it}^* =\{(e_{i1}^*,e_{i2}^*,\cdots,e_{it_i}^*)(\alpha
e_{i1}^*,\alpha e_{it_i}^*,\cdots, e_{i2}^*)\}$ and
$e_{ij}^*=e_{pq}$. Take $e_{ij}^*=e_{pq}$ if  $ q\equiv 1(mod 2)$
and $e_{ij}^*=\alpha e_{pq}$ if $q\equiv 0(mod 2)$.  Then we get
that $M^{\xi\alpha}= M$.}

Now if $o(\xi^*)\equiv 1(mod 2)$, by Lemma $2.7$,
$o(\xi^*\alpha)=2o(\xi^*)$. Therefore, for a chosen quadricell in
$(e^{i1},e^{i2},\cdots,e^{it_i})$ adjacent to the vertex $x_{i1}$
for $i=1,2,\cdots,n$, where, $n=$ the order of the graph $\Gamma$,
the resultant map $M$ is unstable under the action of $\xi\alpha$.
Whence, $\xi\alpha$ is not an automorphism of a map underlying
$\Gamma$. \quad\quad $\natural$

\vskip 8mm

\no{\bf $3$. Determine automorphisms of complete maps} \vskip 6mm

\no Now we determine all automorphisms of complete maps in this
section by applying the results gotten in Section $2$.

Let $K_n$ be a complete graph of order $n$. Label its vertices by
integers $1,2,...,n$. Then its edge set is $\{ij| 1\leq i,j\leq n
,i\not=j\quad and\quad ij=ji\}$. For convenience, we use $i^j$
denoting an edge $ij$ of the complete graph $K_n$ and $i^j = j^i ,
1\leq i,j\leq n , i\not=j$. Then its quadricells of this edge can
be represented by $\{ i^{j+} , i^{j-} , j^{i+ }, j^{i-} \}$ and

$${\mathcal X}_{\alpha,\beta} (K_n)=\{ i^{j+} : 1\leq i,j\leq n ,i\not=j \}
\bigcup\{ i^{j-} : 1\leq i,j\leq n ,i\not=j \},$$

$$\alpha = \prod\limits_{1\leq i,j\leq n ,i\not=j}  ( i^{j+} , i^{j-}),$$

$$\beta = \prod\limits_{1\leq i,j\leq n ,i\not=j}  ( i^{j+} , i^{j+})( i^{j-} , i^{j-}).$$

Recall that the automorphism group of $K_n$ is just the symmetry
group of degree $n$, i.e., ${\rm AutK}_n=S_n$. The above
representation enables us to determine all automorphisms of
complete maps of order $n$ on surfaces.

\vskip 4mm \no{\bf Theorem $3.1$} {\it All cyclic order-preserving
automorphisms of non-orientable complete maps of order$\geq 4$ are
extended actions of elements in}

$${\mathcal E}_{[s^{\frac{n}{s}}]}, \quad {\mathcal E}_{[1,s^{\frac{n-1}{s}}]},$$

\no{and all cyclic order-reversing automorphisms of non-orientable
complete maps of order$\geq 4$ are extended actions of elements
in}

$$
\alpha {\mathcal E}_{[(2s)^{\frac{n}{2s}}]},\quad \alpha {\mathcal
E}_{[(2s)^{\frac{4}{2s}}]}, \quad\alpha {\mathcal E}_{[1,1,2]},
$$
\no{where, ${\mathcal E}_{\theta}$ denotes the conjugate class
containing element $\theta$ in the symmetry group $S_n$} \vskip
3mm

{\it Proof} \ Firstly, we prove that the induced permutation
$\xi^*$ on complete map of order $n$ by an element $\xi\in S_n$ is
an cyclic order-preserving automorphism of a non-orientable map,
if, and only if,

$$\xi\in {\mathcal E}_{s^{\frac{n}{s}}}\bigcup {\mathcal E}_{[1,s^{\frac{n-1}{s}}]}$$

Assume the cycle index of $\xi$ is
$[1^{k_1},2^{k_2},...,n^{k_n}]$. If there exist two integers
$k_i,k_j \not=0$, and $i,j\geq 2, i\not=j$, then in the cycle
decomposition of $\xi$ , there are two cycles

$$(u_1,u_2,...,u_i)\quad  {\rm and}\quad  (v_1,v_2,...,v_j).$$

Since
$$[l^{\xi}(u_1), l^{\xi}(u_2)]=i\quad  {\rm and}\quad  [l^{\xi}(v_1), l^{\xi}(v_2)]=j$$

and $i\not= j$, we know that $\xi^*$ is not an automorphism of
embedding by Lemma $2.5$. Whence, the cycle index of $\xi$ must be
the form of $[1^k,s^l]$.

Now if $k\geq 2$, let $(u),(v)$ be two cycles of length $1$ in the
cycle decomposition of $\xi$. By Lemma $2.5$, we know that

$$[l^{\xi}(u), l^{\xi}(v)] = 1.$$

If there is a cycle $(w,...)$ in the cycle decomposition of $\xi$
whose length greater or equal to two, we get that

$$[l^{\xi}(u), l^{\xi}(w)]=[1, l^{\xi}(w)]= l^{\xi}(w).$$

According to Lemma $2.5$, we get that $l^{\xi}(w)=1$, a
contradiction. Therefore, the cycle index of $\xi$ must be the
forms of $[s^l]$ or $[1, s^l]$. Whence, $ sl=n$ or $sl+1=n$.
Calculation shows that $l=\frac{n}{s}$ or $l=\frac{n-1}{s}$. That
is, the cycle index of $\xi$ is one of the following three types
$[1^n]$,  $[1,s^{\frac{n-1}{s}}]$  and   $[s^{\frac{n}{s}}]$  for
some integer $s$ .

Now we only need to prove that for each element $\xi$ in
${\mathcal E}_{[1,s^{\frac{n-1}{s}}]}$  and ${\mathcal
E}_{[s^{\frac{n}{s}}]}$ , there exists an non-orientable complete
map $M$ of order $n$ with an induced permutation $\xi^*$ being its
cyclic order-preserving automorphism of surface. The discussion
are divided into two cases.

\vskip 2mm {\bf Case 1}   \hspace{40mm}    $ \xi\in{ {\mathcal
E}_{[s^{\frac{n}{s}}]}}$ \vskip 2mm Assume the cycle decomposition
of $\xi$ being $\xi =
(a,b,\cdots,c)\cdots(x,y,\cdots,z)\cdots(u,v,$ $\cdots,w)$, where,
the length of each cycle is $k$, and $1\leq
a,b,\cdots,c,x,y,\cdots,z,u,v,\cdots,w\leq n$ . In this case, we
can construct a non-orientable complete map $M_{1} = ({\mathcal
X}_{\alpha,\beta}^{1},{\mathcal{P}}_{1})$ as follows.

$${\mathcal X}_{\alpha,\beta}^{1}=\{ i^{j+} : 1\leq i,j\leq n ,i(j \}\bigcup
\{ i^{j-} : 1\leq i,j\leq n ,i\not=j \},$$

$${\mathcal{P}}_{1}= \prod\limits_{x\in\{a,b,\cdots,c,\cdots,x,y,\cdots,z,u,v,\cdots,w\}}
(C(x))(\alpha C(x)^{-1}\alpha),$$

\no{where,}

$$C(x) = (x^{a+},\cdots ,x^{x*},\cdots,x^{u+},x^{b+},x^{y+},\cdots,
\cdots,x^{v+},x^{c+},\cdots,x^{z+},\cdots ,x^{w+}),$$

\no{$x^{x*}$ denotes an empty position and }

$$\alpha C(x)^{-1}\alpha =(x^{a-},x^{w-},\cdots,x^{z-},\cdots ,x^{c-},x^{v-},\cdots,
x^{b-},x^{u-},\cdots,x^{y-},\cdots ).$$

It is clear that $M_1^{\xi^*} = M_1$. Therefore, $\xi^*$ is an
cyclic order-preserving automorphism of the map $M_1$.

\vskip 3mm {\bf Case 2}  \hspace{40mm}  $\xi \in{ {\mathcal
E}_{[1,s^{\frac{n-1}{s}}]}}$ \vskip 3mm We assume the cycle
decomposition of $\xi$ being

$$\xi = (a,b,...,c)...(x,y,...,z)...(u,v,...,w)(t),$$

\no{where, the length of each cycle is $k$ beside the final cycle,
and $1\leq a,b...c, x,y..., z,$ $u,v,...,w,t\leq n$ . In this
case, we construct a non-orientable complete map $M_{2} =
({\mathcal X}_{\alpha,\beta}^{2},{\mathcal{P}}_{2})$
 as follows.}

$${\mathcal X}_{\alpha,\beta}^{2}=\{ i^{j+} : 1\leq i,j\leq n ,i\not=j \}\bigcup
\{ i^{j-} : 1\leq i,j\leq n ,i\not=j \},$$

$${\mathcal{P}}_{2} = (A)( \alpha A^{-1})\prod\limits_{x\in\{a,b,...,c,...,
x,y,...z,u,v,...,w\}} (C(x))(\alpha C(x)^{-1}\alpha),$$

\no{where,}

$$A= (t^{a+},t^{x+},...t^{u+},t^{b+},t^{y+},...,t^{v+},...,
t^{c+},t^{z+},...,t^{w+}),$$

$$\alpha A^{-1}\alpha=(t^{a-},t^{w-},...t^{z-},t^{c-},t^{v-},...,t^{y-},...,
t^{b-},t^{u-},...,t^{x-}),$$

$$C(x) = (x^{a+},...,x^{x*},...,x^{u+},x^{b+},...,x^{y+},...,x^{v+},...,x^{c+},...,
x^{z+},...,x^{w+})$$

\no{and}

$$\alpha C(x)^{-1}\alpha =(x^{a-},x^{w-},..,x^{z-},...,x^{c-},...,x^{v-},...,
x^{y-},...,x^{b-},x^{u-},...).$$

It is also clear that $M_2^{\xi^*} = M_2$. Therefore, $\xi^*$ is
an automorphism of the map $M_2$ . \vskip 2mm

Now we consider the case of cyclic order-reversing automorphisms
of a complete map. According to Lemma $2.8$, we know that an
element $\xi\alpha$, where, $\xi\in S_n$, is an cyclic
order-reversing automorphism of a complete map only if,

$$\xi\in {\mathcal E}_{[k^{\frac{n_1}{k}},(2k)^{\frac{n-n_1}{2k}}]}.$$

Our discussion is divided into two parts. \vskip 3mm \no{\bf Case
$3$}\quad\quad\quad\quad   $n_1 = n$ \vskip 3mm

Without loss of generality, we can assume the cycle decomposition
of $\xi$ has the following form in this case.

$$
\xi
=(1,2,\cdots,k)(k+1,k+2,\cdots,2k)\cdots(n-k+1,n-k+2,\cdots,n).
$$

\vskip 3mm \no{\bf Subcase $3.1$}\quad\quad\quad $k\equiv 1(mod
2)$ and $k > 1$ \vskip 3mm

According to Lemma $2.8$, we know that $\xi^*\alpha$ is not an
automorphism of maps since $o(\xi^*)=k\equiv 1(mod 2)$.

\vskip 3mm \no{\bf Subcase $3.2$} \quad\quad\quad\quad  $k\equiv
0(mod 2)$ \vskip 3mm Construct a non-orientable map
$M_3=({\mathcal X}_{\alpha,\beta}^3,{\mathcal P}_3)$, where $ X^3=
E(K_n)$ and

$$
{\mathcal P}_3=\prod\limits_{i\in \{1,2,\cdots,n\}}(C(i))(\alpha
C(i)^{-1}\alpha),
$$

\no{where, if $i\equiv 1(mod 2)$, then}

$$C(i)=(i^{1+},i^{k+1+},\cdots,i^{n-k+1+},i^{2+},\cdots,i^{n-k+2+},\cdots,i^{i*},\cdots,i^{k+},i^{2k+},\cdots,i^{n+}),$$

$$\alpha C(i)^{-1}\alpha=(i^{1-},i^{n-},\cdots,i^{2k-},i^{k-},\cdots,i^{k+1-})$$

\no{and if $i\equiv 0(mod 2)$, then}

$$C(i)=(i^{1-},i^{k+1-},\cdots,i^{n-k+1-},i^{2-},\cdots,i^{n-k+2-},\cdots,i^{i*},\cdots,i^{k-},i^{2k-},\cdots,i^{n-}),$$

$$\alpha C(i)^{-1}\alpha=(i^{1+},i^{n+},\cdots,i^{2k+},i^{k+},\cdots,i^{k+1+}).$$

Where, $i^{i*}$ denotes the empty position, for example,
$(2^1,2^{2*},2^3,2^4,2^5)=(2^1,2^3,2^4,2^5)$. It is clear that
${\mathcal P}_3^{\xi\alpha}={\mathcal P}_3$, that is, $\xi\alpha$
is an automorphism of map $M_3$.

\vskip 3mm \no{\bf Case 4}\quad\quad\quad\quad $n_1\not= n$ \vskip
3mm

Without loss of generality, we can assume that

\begin{eqnarray*}
\xi &=& (1,2,\cdots,k)(k+1,k+2,\cdots,n_1)\cdots (n_1-k+1,n_1-k+2,\cdots,n_1)\\
&\times&
(n_1+1,n_1+2,\cdots,n_1+2k)(n_1+2k+1,\cdots,n_1+4k)\cdots(n-2k+1,\cdots,n)
\end{eqnarray*}

\vskip 3mm \no{\bf Subcase $4.1$} \quad\quad\quad\quad $k\equiv
0(mod 2)$ \vskip 3mm

Consider the orbits of $1^{2+}$ and $n_1+2k+1^{1+}$ under the
action of $<\xi\alpha >$, we get that

$$|orb((1^{2+})^{<\xi\alpha >})| = k$$

\no{and}

$$|orb(((n_1+2k+1)^{1+})^{< \xi\alpha >})| = 2k.$$

\no{Contradicts to Lemma $2.5$.}

\vskip 3mm \no{\bf Subcase $4.2$} \quad\quad\quad\quad $k\equiv
1(mod 2)$ \vskip 3mm

In this case, if $k\not= 1$, then $k\geq 3$. Similar to the
discussion of Subcase $3.1$, we know that $\xi\alpha$ is not an
automorphism of complete map. Whence, $k=1$ and

$$\xi\in {\mathcal E}_{[1^{n_1},2^{n_2}]}.$$

\no{Without loss of generality, assume that}

$$\xi =(1)(2)\cdots(n_1)(n_1+1,n_1+2)(n_1+3,n_1+4)\cdots(n_1+n_2-1,n_1+n_2).$$

If $n_2\geq 2$, and there exists a map $M=({\mathcal
X}_{\alpha,\beta},{\mathcal P})$, assume the vertex $v_1$ in $M$
being

$$
v_1=(1^{l_{12}+},1^{l_{13}+},\cdots,1^{l_{1n}+})(1^{l_{12}-},1^{l_{1n}-},\cdots,1^{l_{13}-})
$$
\no{where, $l_{1i}\in\{+2,-2,+3,-3,\cdots,+n,-n\}$ and
$l_{1i}\not= l_{1j}$ if $i\not= j$.}

Then we get that

$$
(v_1)^{\xi\alpha}=(1^{l_{12}-},1^{l_{13}-},\cdots,1^{l_{1n}-})
(1^{l_{12}+},1^{l_{1n}+},\cdots,1^{l_{13}+})\not=v_1.
$$

\no{Whence, $\xi\alpha$ is not an automorphism of map $M$, a
contradiction.

Therefore, $n_2 =1$. Similarly, we can also get that $n_1 =2$.
Whence,  $\xi =(1)(2)(34)$ and $n=4$. We construct a stable
non-orientable map $M_4$ under the action of $\xi\alpha$ as
follows.

$$M_4=({\mathcal X}_{\alpha,\beta}^4,{\mathcal P}_4),$$

\no{where,}

\begin{eqnarray*}
{\mathcal P}_4  &=& (1^{2+},1^{3+},1^{4+})(2^{1+},2^{3+},2^{4+})(3^{1+},3^{2+},3^{4+})(4^{1+},4^{2+},4^{3+})\\
& \times &
(1^{2-},1^{4-},1^{3-})(2^{1-},2^{4-},2^{3-})(3^{1-},3^{4-},3^{2-})(4^{1-},4^{3-},4^{2-}).
\end{eqnarray*}

Therefore, all cyclic order-preserving automorphisms of
non-orientable complete maps are extended actions of elements in

$${\mathcal E}_{[s^{\frac{n}{s}}]}, \quad {\mathcal E}_{[1,s^{\frac{n-1}{s}}]}$$

\no{and all cyclic order-reversing automorphisms of non-orientable
complete maps are extended actions of elements in}

$$
 \alpha {\mathcal E}_{[(2s)^{\frac{n}{2s}}]},\quad \alpha {\mathcal E}_{[(2s)^{\frac{4}{2s}}]}
\quad  \alpha {\mathcal E}_{[1,1,2]}.
$$
\no{This completes the proof.\quad\quad $\natural$}

\vskip 3mm According to the Rotation Embedding Scheme for
orientable embedding of a graph formalized by Edmonds in $[5]$,
each orientable complete map is just the case of eliminating the
signs "+, -" in our representation for complete maps. Whence,we
also get the following result for automorphisms of orientable
complete maps, which is similar to Theorem $3.1$.

\vskip 4mm \no{\bf Theorem $3.2$} {\it All cyclic order-preserving
automorphisms of orientable complete maps of order$\geq 4$ are
extended actions of elements in}

$${\mathcal E}_{[s^{\frac{n}{s}}]}, \quad {\mathcal E}_{[1,s^{\frac{n-1}{s}}]} $$

\no{ and all cyclic order-reversing automorphisms of orientable
complete maps of order$\geq 4$ are extended actions of elements in
}
$$
\alpha {\mathcal E}_{[(2s)^{\frac{n}{2s}}]},\quad \alpha {\mathcal
E}_{[(2s)^{\frac{4}{2s}}]}, \quad\alpha {\mathcal E}_{[1,1,2]},
$$
\no{where,${\mathcal E}_{\theta}$ denotes the conjugate class
containing $\theta$ in $S_n$.} \vskip 3mm

{\it Proof} \ The proof is similar to that of Theorem $3.1$. For
completion, we only need to construct orientable maps $M_i^O ,
i=1,2,3,4$ to replace these non-orientable maps $M_1, i=1,2,3,4$
in the proof of Theorem $3.1$.

In fact, for cyclic order-preserving case, we only need to take
$M_1^O$, $M_2^O$ to be the resultant maps eliminating the signs +
- in $M_1$, $M_2$ constructed in the proof of Theorem $3.1$.

For  the cyclic order-reversing case, we take  $M_3^O=
(E(K_n)_{\alpha ,\beta},{\mathcal P}_3^O)$ with

$$
{\mathcal P}_3=\prod\limits_{i\in \{1,2,\cdots,n\}}(C(i)),
$$

\no{where, if $i\equiv 1(mod 2)$, then}

$$C(i)=(i^{1},i^{k+1},\cdots,i^{n-k+1},i^{2},\cdots,i^{n-k+2},\cdots,i^{i*},\cdots,i^{k},i^{2k},\cdots,i^{n}),$$

\no{and if $i\equiv 0(mod 2)$, then}

$$C(i)=(i^{1},i^{k+1},\cdots,i^{n-k+1},i^{2},\cdots,i^{n-k+2},\cdots,i^{i*},\cdots,i^{k},i^{2k},\cdots,i^{n})^{-1},$$

\no{where $i^{i*}$ denotes the empty position and $M_4^O
=(E(K_4)_{\alpha ,\beta}, {\mathcal P}_4)$ with}

$$
{\mathcal P}_4=
(1^2,1^3,1^4)(2^1,2^3,2^4)(3^1,3^4,3^2)(4^1,4^2,4^3).
$$

It can be shown that $(M_i^O)^{g*}=M_i^O$, $i=1,2$ and
$(M_i^O)^{\xi\alpha}=M_i^O$ for $i=3,4$. $\natural$

All results in this section are useful for the enumeration of
complete maps in the next section.

\vskip 10mm

\no{\bf 4. The Enumeration of complete maps on surfaces} \vskip
5mm

\no We first consider the permutation and its stabilizer . The
permutation with the following form $(x_1,x_2,...,x_n)( \alpha
x_n, \alpha x_2,...,\alpha x_1)$ is called a pair permutation. The
following result is obvious.

\vskip 4mm \no{{\bf Lemma 4.1} {\it Let $g$ be a permutation on
the set $\Omega=\{x_1,x_2,...,x_n\}$ such that  $g\alpha=\alpha g
$. If}}

$$g(x_1,x_2,...,x_n)( \alpha x_n, \alpha x_{n-1},..., \alpha x_1)g^{-1}=
(x_1,x_2,...,x_n)( \alpha x_n, \alpha x_{n-1},..., \alpha x_1),$$

\no{\it then}

$$
g=(x_1,x_2,...,x_n)^k
$$

\no{\it and if}

$$g\alpha (x_1,x_2,...,x_n)( \alpha x_n, \alpha x_{n-1},..., \alpha x_1)(g\alpha )^{-1}=
(x_1,x_2,...,x_n)( \alpha x_n, \alpha x_{n-1},..., \alpha x_1),$$

\no{\it then}

$$
g\alpha =(\alpha x_n, \alpha x_{n-1},..., \alpha x_1)^k
$$

\no{for some integer $k,1\leq k\leq n$.}

\vskip 4mm \no{{\bf Lemma 4.2} {\it For each permutation $g,g\in
{\mathcal E}_{[k^{\frac{n}{k}}]}$ satisfying $g\alpha=\alpha g$ on
the set $\Omega=\{x_1,x_2,...,x_n\}$, the number of stable pair
permutations in $\Omega$ under the action of $g$ or $g\alpha$ is}}

$$
\frac{2\phi (k)(n-1)!}{|{\mathcal E}_{[k^{\frac{n}{k}}]}|},
$$

\no{\it  where  $\phi (k)$ denotes the Euler function.} \vskip 3mm

{\it Proof} \ Denote the number of stable pair permutations under
the action of $g$ or $g\alpha$ by $n(g)$ and $\mathcal C$ the set
of pair permutations. Define the set $A=\{(g,C)| g\in {\mathcal
E}_{[k^{\frac{n}{k}}]}, C\in{\mathcal C}\quad {\rm and} \quad
C^g=C \ or \ $ $C^{g\alpha}=C\}$. Clearly, for $\forall
g_1,g_2\in{ {\mathcal E}_{[k^{\frac{n}{k}}]}}$, we have
$n(g_1)=n(g_2)$. Whence, we get that

$$
|A| = | {\mathcal E}_{[k^{\frac{n}{k}}]}| n(g) . \hspace{30mm}
(4.1)
$$

On the other hand, by Lemma 4.1, for any pair permutation
$C=(x_1,x_2,...,x_n)$ $( \alpha x_n, \alpha x_{n-1},..., \alpha
x_1)$, since $C$ is stable under the action of $g$, there must be
$g=(x_1,x_2,...,x_n)^l $ or $g\alpha =(\alpha x_n, \alpha
x_{n-1},..., \alpha x_1)^l$, where $l=s\frac{n}{k}, 1\leq s\leq k$
and $(s,k)=1$. Therefore, there are $2\phi(k)$ permutations in
${\mathcal E}_{[k^{\frac{n}{k}}]}$ acting on it stable. Whence, we
also have

$$
|A| = 2\phi (k) |{\mathcal C}|. \hspace{30mm}         (4.2)
$$

Combining (4.1) with (4.2), we get that

$$
n(g)=\frac{2\phi (k)|{\mathcal C}|}{|{\mathcal
E}_{[k^{\frac{n}{k}}]}|}= \frac{2\phi (k)(n-1)!}{ |{\mathcal
E}_{[k^{\frac{n}{k}}]}| }.\hspace{20mm}\natural
$$

\vskip 3mm

 Now we can enumerate the unrooted complete maps on surfaces.

\vskip 4mm \no{{\bf Theorem 4.1} {\it The number $n^{L}(K_n)$ of
complete maps of order $n\geq 5$ on surfaces is}}

$$
n^{L}(K_n)=\frac{1}{2}(\sum\limits_{k|n}+\sum\limits_{k|n, k\equiv
0(mod 2)})
 \frac{2^{\alpha (n,k)}(n-2)!^{\frac{n}{k}}}
{k^{\frac{n}{k}}(\frac{n}{k})!}+ \sum\limits_{k|(n-1),k\not=1}
\frac{\phi (k)2^{\beta (n,k)}(n-2)!^{\frac{n-1}{k}}}{n-1},
$$
\no{where, }

\[
\alpha (n,k)=\left\{\begin{array}{cc}
\frac{n(n-3)}{2k}, & {\rm if}\quad k\equiv 1(mod 2);\\
\frac{n(n-2)}{2k}, & {\rm if}\quad k\equiv 0(mod 2),
\end{array}
\right.
\]

\no{and}

\[
\beta (n,k)=\left\{\begin{array}{cc}
\frac{(n-1)(n-2)}{2k}, & {\rm if}\quad k\equiv 1(mod 2);\\
\frac{(n-1)(n-3)}{2k}, & {\rm if}\quad k\equiv 0(mod 2).
\end{array}
\right.
\]
\no{and $n^L(K_4)=11.$}

\vskip 5mm

{\it Proof} \  According to ($2.3$) in Corollary $2.1$ and Theorem
$3.1$ for $n\geq 5$, we know that

\begin{eqnarray*}
n^{L}(K_n) &=& \frac{1}{2|{\rm AutK}_n|}
\times(\sum\limits_{g_1\in{{\mathcal E}_{[k^{\frac{n}{k}}]}}}|\Phi
(g_1)| + \sum\limits_{g_2\in{{\mathcal
E}_{[(2s)^{\frac{n}{2s}}]}}}|\Phi (g_2\alpha)|\\
&+& \sum\limits_{h\in{{\mathcal E}_{[1,k^{\frac{n-1}{k}}]}}}|\Phi (h)|) \\
&=& \frac{1}{2n!}\times (\sum\limits_{k|n}|{\mathcal
E}_{[k^{\frac{n}{k}}]}| |\Phi (g_1)|+ \sum\limits_{l|n,l\equiv
0(mod2)}|{\mathcal E}_{[l^{\frac{n}{l}}]}| |\Phi (g_2\alpha)|\\
&+& \sum\limits_{l|(n-1)}|{\mathcal E}_{[1,l^{\frac{n-1}{l}}]}|
|\Phi (h)|),
\end{eqnarray*}

\no{where, $g_1\in{{\mathcal
E}_{[k^{\frac{n}{k}}]}},g_2\in{{\mathcal E}_{[l^{\frac{n}{l}}]}}$
and $h\in{{\mathcal E}_{[1,k^{\frac{n-1}{k}}]}}$ are three chosen
elements. }

Without loss of generality, we assume that an element
$g,g\in{{\mathcal E}_{[k^{\frac{n}{k}}]}}$ has the following cycle
decomposition.

$$g = (1,2,...,k)(k+1,k+2,...,2k)...((\frac{n}{k}-1)k+1,(\frac{n}{k}-1)k+2,...,n)$$

and

$$
{\mathcal P}={\prod}_1\times{\prod}_2,
$$

where,

$$
{\prod}_1=(1^{i_{21}},1^{i_{31}},...,1^{i_{n1}})(2^{i_{12}},2^{i_{32}},...,2^{i_{n2}})...
(n^{i_{1n}},n^{i_{2n}},...,n^{i_{(n-1)n}}),
$$

and

$${\prod}_2=\alpha ({{\prod}_1}^{-1})\alpha^{-1}$$

\no{being a complete map which is stable under the action of $g$, where
$s_{ij}\in \{k+,k-| k= 1,2,...n \}$.}

Notice that the quadricells adjacent to the vertex "$1$" can make
$2^{n-2}(n-2)!$ different pair permutations and for each chosen
pair permutation, the pair permutations adjacent to the vertices
$2,3,...,k$ are uniquely determined since $\mathcal P$  is stable
under the action of $g$.

Similarly, for each given pair permutation adjacent to the vertex
$k+1,2k+1,..., (\frac{n}{k}-1)k+1$, the pair permutations adjacent
to $k+2,k+3,...,2k$ and $2k+2,2k+3,...,3k$ and,...,and
$(\frac{n}{k}-1)k+2,(\frac{n}{k}-1)k+3,...n$ are also uniquely
determined because $\mathcal P$ is stable under the action of $g$.

Now for an orientable embedding $M_1$ of $K_n$, all the induced
embeddings by exchanging two sides of some edges and retaining the
others unchanged in $M_1$ are the same as $M_1$ by the definition
of maps. Whence, the number of different stable embeddings under
the action of $g$ gotten by exchanging $x$ and $\alpha x$ in $M_1$
for $x\in U, U\subset{\mathcal X}_{\beta}$, where ${\mathcal
X}_{\beta}= \bigcup\limits_{x\in E(K_n)} \{x,\beta x \}$ , is
$2^{g(\varepsilon)-\frac{n}{k}}$, where $g(\varepsilon)$ is the
number of orbits of $E(K_n)$ under the action of $g$ and we
substract $\frac{n}{k}$ because we can chosen
$1^{2+},k+1^{1+},2k+1^{1+},\cdots,n-k+1^{1+}$ first in our
enumeration.

Notice that the length of each orbit under the action of $g$ is $k$ for
$\forall x\in E(K_n)$ if $k$ is odd and is $\frac{k}{2}$ for
$x=i^{i+\frac{k}{2}}, i=1,k+1,\cdots, n-k+1$, or $k$ for all other edges
if $k$ is even. Therefore, we get that

\[
g(\varepsilon)=\left\{\begin{array}{cc}
\frac{\varepsilon (K_n)}{k}, & {\rm if}\quad k\equiv 1(mod 2);\\
\frac{\varepsilon (K_n)-\frac{n}{2}}{k}, & {\rm if}\quad k\equiv 0(mod 2).
\end{array}
\right.
\]

\no{Whence, we have that}

\[
\alpha (n,k)=g(\varepsilon) - \frac{n}{k}=\left\{\begin{array}{cc}
\frac{n(n-3)}{2k}, & {\rm if}\quad k\equiv 1(mod 2);\\
\frac{n(n-2)}{2k}, & {\rm if}\quad k\equiv 0(mod 2),
\end{array}
\right.
\]

\no{and}

$$
|\Phi (g)|=2^{\alpha (n,k)}(n-2)!^{\frac{n}{k}},\hspace{30mm} (4.3)
$$

Similarly, if $k\equiv 0(mod2)$, we get also that

$$
|\Phi (g\alpha)|=2^{\alpha (n,k)}(n-2)!^{\frac{n}{k}}
\hspace{30mm} (4.4)
$$

\no{for an chosen element $g$, $g\in {\mathcal
E}_{[k^{\frac{n}{k}}]}.$}

Now for $\forall h\in{{\mathcal E}_{[1,k^{\frac{n-1}{k}}]}} $,
without loss of generality, we assume that $h =
(1,2,...,k)(k+1,k+2,...,2k)...
((\frac{n-1}{k}-1)k+1,(\frac{n-1}{k}-1)k+2,...,(n-1))(n)$. Then
the above statement is also true for the complete graph $K_{n-1}$
with the vertices $1,2,\cdots,n-1$. Notice that the quadricells
$n^{1+},n^{2+},\cdots,n^{n-1+}$ can be chosen first in our
enumeration and they are not belong to the graph $K_{n-1}$.
According to Lemma 4.2, we get that

$$
|\Phi (h)|=2^{\beta (n,k)}(n-2)!^{\frac{n-1}{k}}\times\frac{2\phi
(k)(n-2)!} {| {\mathcal E}_{[1,k^{\frac{n-1}{k}}]}|},\hspace{10mm}
(4.5)
$$

\no{Where}

\[
\beta (n,k)=h(\varepsilon)=\left\{\begin{array}{cc}
\frac{\varepsilon (K_{n-1})}{k}-\frac{n-1}{k}=\frac{(n-1)(n-4)}{2k},
& {\rm if}\quad k\equiv 1(mod 2);\\
\frac{\varepsilon (K_{n-1})}{k}-\frac{n-1}{k}=\frac{(n-1)(n-3)}{2k},
& {\rm if}\quad k\equiv 0(mod 2).
\end{array}
\right.
\]

Combining $(4.3)-(4.5)$, we get that

\begin{eqnarray*}
n^{L}(K_n) &=&  \frac{1}{2n!}\times (\sum\limits_{k|n}|{\mathcal
E}_{[k^{\frac{n}{k}}]}| |\Phi (g_0)|+\sum\limits_{l|n,l\equiv
0(mod2)}|{\mathcal E}_{[l^{\frac{n}{l}}]}| |\Phi (g_1\alpha)|\\
&+& \sum\limits_{l|(n-1)}|{\mathcal E}_{[1,l^{\frac{n-1}{l}}]}|
|\Phi
(h)|)\\
 &=& \frac{1}{2n!}\times ( \sum\limits_{k|n} \frac{n! 2^{\alpha
(n,k)}(n-2)!^{\frac{n}{k}}}{k^{\frac{n}{k}}(\frac{n}{k})!}+
\sum\limits_{k|n,k\equiv 0(mod2)}
\frac{n! 2^{\alpha (n,k)}(n-2)!^{\frac{n}{k}}}{k^{\frac{n}{k}}(\frac{n}{k})!}\\
&+&\sum\limits_{k|(n-1),k\not =1}
\frac{n!}{k^{\frac{n-1}{k}}(\frac{n-1}{k})!} \times\frac{2\phi
(k)(n-2)! 2^{\beta (n,k)}(n-2)!^{\frac{n-1}{k}}}
{\frac{(n-1)!}{k^{\frac{n-1}{k}}(\frac{n-1}{k})!}})\\
&=& \frac{1}{2}(\sum\limits_{k|n}+\sum\limits_{k|n,k\equiv
0(mod2)}) \frac{2^{\alpha (n,k)}(n-2)!^{\frac{n}{k}}}
{k^{\frac{n}{k}}(\frac{n}{k})!}+ \sum\limits_{k|(n-1),k\not=1}
\frac{\phi (k)2^{\beta (n,k)}(n-2)!^{\frac{n-1}{k}}}{n-1}.
\end{eqnarray*}

\vskip 3mm

 For $n=4$, similar calculation shows that $n^L(K_4)=11$
by consider the fixing set of permutations in ${\mathcal
E}_{[s^{\frac{4}{s}}]}$,${\mathcal E}_{[1,s^{\frac{3}{s}}]}$,
${\mathcal E}_{[(2s)^{\frac{4}{2s}}]}$,$\alpha{\mathcal
E}_{[(2s)^{\frac{4}{2s}}]}$ and $\alpha {\mathcal E}_{[1,1,2]}$.
\quad\quad $\natural$

 \vskip 3mm

For orientable complete maps, we get the number $n^{O}(K_n)$ of
orientable complete maps of order $n$ as follows.

\vskip 4mm \no{{\bf Theorem 4.2} {\it The number $n^O((K_n)$ of
complete maps of order $n\geq 5$ on orientable surfaces is}}

$$
n^{O}(K_n)=\frac{1}{2}(\sum\limits_{k|n}+\sum\limits_{k|n, k\equiv
0(mod 2)} ) \frac{(n-2)!^{\frac{n}{k}}}
{k^{\frac{n}{k}}(\frac{n}{k})!}+ \sum\limits_{k|(n-1),k\not=1}
\frac{\phi (k)(n-2)!^{\frac{n-1}{k}}}{n-1} .
$$

\no{and $n(K_4)=3.$}

\vskip 3mm {\it Proof}  According to the Tutte's algebraic
representation of maps, a map $M=({\mathcal
X}_{\alpha,\beta},{\mathcal P})$ is orientable if and only if for
$\forall x\in {\mathcal X}_{\alpha ,\beta}$, $x$ and $\alpha\beta
x$ are in a same orbit of ${\mathcal X}_{\alpha ,\beta}$ under the
action of the group $\Psi_I=\left<\alpha\beta,{\mathcal
P}\right>$. Now applying $(2.1)$ in Corollary $2.1$ and Theorem
$3.1$, similar to the proof of Theorem $4.1$, we get the number
$n^O(K_n)$ for $n\geq 5$ as follows}

$$
n^{O}(K_n)=\frac{1}{2}(\sum\limits_{k|n}+\sum\limits_{k|n, k\equiv
0(mod 2)} ) \frac{(n-2)!^{\frac{n}{k}}}
{k^{\frac{n}{k}}(\frac{n}{k})!}+ \sum\limits_{k|(n-1),k\not=1}
\frac{\phi (k)(n-2)!^{\frac{n-1}{k}}}{n-1} .
$$

\no{and for the complete graph $K_4$, calculation shows that
$n(K_4)=3. \quad\quad \natural$}

\vskip 3mm

Notice that $n^O(K_n)+n^N(K_n)=n^L(K_n)$. Therefore, we can also
get the number $n^N(K_n)$ of unrooted complete maps of order $n$
on non-orientable surfaces by Theorem $4.1$ and Theorem $4.2$.

\no{{\bf Theorem 4.3} {\it The number $n^N(K_n)$ of unrooted
complete maps of order $n, n\geq 5$ on non-orientable surfaces
is}}

\begin{eqnarray*}
n^{N}(K_n)&=& \frac{1}{2}(\sum\limits_{k|n}+\sum\limits_{k|n,
k\equiv 0(mod 2)} ) \frac{(2^{\alpha
(n,k)}-1)(n-2)!^{\frac{n}{k}}}
{k^{\frac{n}{k}}(\frac{n}{k})!}\\
&+& \sum\limits_{k|(n-1),k\not=1} \frac{\phi (k)(2^{\beta
(n,k)}-1)(n-2)!^{\frac{n-1}{k}}}{n-1},
\end{eqnarray*}

\no{and $n^N(K_4)=8$. Where, $\alpha (n,k)$ and $\beta(n,k)$ are
same as in Theorem 4.1.}

For $n=5$, calculation shows that $n^L(K_5)=1080$ and $n^{O}(K_5)
=45$ based on Theorem $4.1$ and $4.2$. For $n=4$, there are 3
unrooted orientable maps and 8 non-orientable maps shown in the
Fig.$2$.

\includegraphics[bb=10 10 200 400]{7mg2.eps}

\C{Fig.$2$} \vskip 3mm \no{All the 11 maps of $K_4$ on surfaces
are non-isomorphic.}

Noticing that for an orientable map $M$, its cyclic
order-preserving automorphisms are just the orientation-preserving
automorphisms of map $M$ by definition. Now consider the action of
cyclic order-preserving automorphisms of complete maps, determined
in Theorem $3.2$ on all orientable embeddings of a complete graph
of order $n$. Similar to the proof of Theorem $4.2$, we can get
the number of non-equivalent embeddings of complete graph of order
$n$, which is same as the result of Mull et al. in $[15]$.

 \vskip 10mm

{\bf References}

\vskip 5mm

\re{[1]}N.L.Biggs, Automorphisms of imbedded graphs, {\it
J.Combinatorical Theory} 11,132-138(1971).

\re{[2]}N.L.Biggs and A.T.White, {\it Permutation Groups and
Combinatoric Structure}, Cambridge University Press (1979).

\re{[3]}E.Bujalance, Automorphism groups of compact planar Klein
surfaces, {\it Manuscripta Math},56,105-124(1986).

\re{[4]}E.Bujalance, Cyclic groups of automorphisms of compact
non-orientable Klein surfaces without boundary, {\it Pacific J.
Math}, vol.109,No.2,279-289(1983).

\re{[5]} Edmonds.J, A combinatorial representation for polyhedral
surfaces, {\it Notices Amer. Math. Soc} 7 (1960)

\re{[6]}W.J.Harvey, Cyclic groups of automorphisms of a compact
Riemann surface, {\it Quart. J.Math.Oxford} (2),17,86-97(1966).

\re{[7]} Jin Ho Kwak and Jaeun Lee, Enumeration of graph
embeddings, {\it Discrete Math}, 135(1994), 129~151.

\re{[8]}L.F. Mao and Y.P. Liu , On the roots on orientable
embeddings of graph, {\it Acta Math. Scienta} Ser.A,vol.23,
No.3,287-293(2003).

\re{[9]}L.F. Mao and Y.P. Liu, Group action for enumerating maps
on surfaces, {\it J.Appl.Math. \&
Computing},vol.13,No.1-2,201-215(2003).

\re{[10]}L.F. Mao and Y.P. Liu, New approach for enumerating maps
on orientable surface, {\it Australasian J.Combinatoric}
(accepted).

\re{[11]}L.F. Mao, Y.P. Liu and F. Tian, Automorphisms of maps
with given underlying graph and its application to enumeration,
{\it Acta Math. Sinica}, vol.21, No.2(2005), 225-236.

\re{[12]} V.A.Liskovets, Enumeration of non-isomorphic planar
maps,{\it Sel.Math.Sov.}, 4;4 (1985), 303~323.

\re{[13]} V.A.Liskovets, A reductive technique for Enumerating
non-isomorphic planar maps,{\it Discrete Math}, 156(1996),
197~217.

\re{[14]} V.A.Liskovets and T.R.S.Walsh, The enumeration of
non-isomorphic 2-connected planar
 maps,{\it Canad.J.Math.}, 35(1983), 417~435.

\re{[15]} B.P.Mull,R.G.Rieper and A.T.White, Enumeration 2-cell
imbeddings of connected graphs,{\it Proc.Amer.Math.Soc.},
103(1988), 321~330.

\re{[16]} B.P.Mull, Enumerating the orientable $2$-cell imbeddings
of complete bipartite graphs, {\it J.Graph Theory}, vol 30,
2(1999),77-90.

\re{[17]} K. Nakagama, On the order of automorphisms of a closed
Riemann surface, {\it Pacific J. Math.},
vol.115,No.2,435-443(1984)

 \re{[18]} S.Negami, Enumeration of
projective -planar embeddings of graphs, {\it Discrete Math},
62(1986), 299~306.

\re{[19]} W.T.Tutte, What is a maps? in {\it New Directions in the
Theory of Graphs} (ed.by F.Harary), Academic Press (1973),
309~325.

\re{[20]} A.T.White, {\it Graphs of Group on Surfaces-
interactions and models}, Elsevier Science B.V. (2001).

\re{[21]}Yanpei Liu, {\it Enumerative Theory of Maps}, Kluwer
Academic Publisher, Dordrecht / Boston / London (1999).

\end{document}